\documentclass{article}
\usepackage{graphicx}
\usepackage{amsmath}
\usepackage{amsfonts}
\usepackage{verbatim}
\usepackage{amssymb}
\usepackage{amsthm}
\usepackage{color}
\definecolor{red}{rgb}{1,0,0}

\definecolor{blue}{rgb}{0,0,1}

\usepackage{url}
\usepackage{hyperref}  
\usepackage[mathlines]{lineno}
\newtheorem{thm}{Theorem}[section]

\newtheorem{ex}[thm]{Example}

\def\ord#1{| #1 |}
\newcommand{\nul}{\operatorname{null}}

\newcommand{\Gc}{\overline{G}}
\newcommand{\G}{\mathcal{G}}
\newcommand{\D}{\Gamma}
\newcommand{\Y}{{\cal Y}}

\newcommand{\R}{\mathbb{R}}
\newcommand{\Z}{\mathbb{Z}}
\newcommand{\C}{\mathbb{C}}
\newcommand{\QQ}{\mathbb{Q}}
\newcommand{\Rnn}{\R^{n\times n}} 
 
\newcommand{\sRnn}{S_n(\R)}

\newcommand{\Fnn}{F^{n\times n}}

\newcommand{\bx}{{\bf x}}
\newcommand{\bv}{{\bf v}}

\newcommand{\bfe}{{\bf e}}

\def\wh#1{\widehat{#1}} 
\newcommand{\SG}{\mathcal{S}(G)}
\newcommand{\SGh}{\mathcal{S}(\wh G)}
\newcommand{\SGp}{\mathcal{S}_+(G)}
\newcommand{\sym}{\mathcal{S}}
\newcommand{\rank}{\operatorname{rank}}
\newcommand{\mr}{\operatorname{mr}}
\newcommand{\mrp}{\operatorname{mr}_+}
\newcommand{\mrs}{\operatorname{mr}_-}
\newcommand{\Mrs}{\operatorname{MR}_-}
\newcommand{\match}{\operatorname{match}}
\newcommand{\mrF}{{\mr^F}}
\newcommand{\M}{\operatorname{M}}
\newcommand{\Mp}{\operatorname{M}_+}

\newcommand{\ZFS}{\operatorname{Z}}
\newcommand{\Zp}{\operatorname{Z}_+}
\newcommand{\Zh}{\wh{\operatorname{Z}}}
\newcommand{\Zl}{\operatorname{Z}_\ell}

\newcommand{\tw}{\operatorname{tw}}

\newcommand{\cc}{\operatorname{cc}}

\newcommand{\ED}{\operatorname{ED}}
\newcommand{\PC}{\operatorname{P}}
\newcommand{\cp}{\mathbin{\scriptscriptstyle\square}}

\newcommand{\amr}{\operatorname{amr}}
\newcommand{\sgn}{{\rm sgn}}
\newcommand{\tri}{\operatorname{tri}}
\newcommand{\sns}{\operatorname{SNS}}

\def\PSDG{{\cal H}^+(G)}

\newcommand{\x}{\times}

\newcommand{\sage}{{\em Sage}}
\newcommand{\bit}{\begin{itemize}}
\newcommand{\eit}{\end{itemize}}
\newcommand{\ben}{\begin{enumerate}}
\newcommand{\een}{\end{enumerate}}
\newcommand{\beq}{\begin{equation}}
\newcommand{\eeq}{\end{equation}}
\newcommand{\bea}{\begin{eqnarray*}}
\newcommand{\eea}{\end{eqnarray*}}
\newcommand{\bpf}{\begin{proof}}
\newcommand{\epf}{\end{proof}\ms}
\newcommand{\ms}{\medskip}
\newcommand{\noi}{\noindent}
\def\al{{\rm mult}}

\newcommand{\h}{\mathcal{H}}
\newcommand{\sr}{\operatorname{sign-rank}}
\newcommand{\uppcc}{\operatorname{upp-cc}}
\newcommand{\Q}{\mathcal{Q}}

\begin{document}

\title{Variants on the minimum rank problem: A survey II\thanks{This 
article is based in part on material prepared for the Banff International 
Research Station workshop,``Theory and Applications of Matrices Described 
by Patterns," and the authors thank  
BIRS for their support.  
}
}
\author{Shaun M. Fallat \thanks{Department of Mathematics and Statistics,
University of Regina, Regina, SK, Canada 
(sfallat@math.uregina.ca). Research supported in part by an
NSERC Discovery Research grant.} \and 
Leslie Hogben\thanks{Department of Mathematics, Iowa State University, 
Ames, IA 50011, USA (lhogben@iastate.edu) and American Institute of 
Mathematics, 360 Portage Ave, 
Palo Alto, CA 94306 (hogben@aimath.org).} }

\maketitle


\begin{abstract}
The  minimum rank problem for a (simple) graph $G$ is to determine   the
smallest possible rank over all real symmetric matrices  whose $ij$th
entry (for $i\neq j$) is nonzero whenever   $\{i,j\}$ is an edge in $G$ and is zero otherwise.  This paper surveys the many developments on the (standard) minimum rank problem and its variants  since the survey paper \cite{FH}.  In particular, positive semidefinite minimum rank, zero forcing parameters, and minimum rank problems for patterns are discussed.
 \end{abstract}



\section{Introduction}\label{sintro} 
Since our survey paper \cite{FH} the volume of work, advances, and
interesting open problems on
many different aspects of the minimum rank of graphs has continued
to expand. Furthermore, since
the 2006 AIM workshop that featured graphs and minimum rank, there
have been 
numerous special sessions, 
minisymposia, and a BIRS workshop emphasizing the topic of minimum rank
of graphs. Consequently, we
felt it was timely to produce an updated survey covering more recent
topics and advances on 
the minimum rank of graphs which is meant to serve as a sequel to the original
survey paper \cite{FH}.

Since this work is follow up reporting, we will not repeat all of the
necessary notation or
terminology that was presented in \cite{FH}, so please consult \cite{FH}
if relevant terms or notation are not spelled out here. However,
we will carefully define key terms and notation used within. 

In general the minimum rank of a graph is simply the smallest rank over
a collection of matrices
that are in some way associated with a given graph $G$. As was outlined
in \cite{FH}, this simple
question has it roots in many different topics in combinatorics and
has been a concern for many
researchers over the years. 
Recently, connections have been found between the related graph parameter zero forcing number and control of quantum systems (see Section \ref{sZ}), and between the minimum rank of sign patterns and communication complexity (see Section \ref{ssSP}). As mentioned above, 
minimum rank problems are a hot  topic  
currently and  
has seen a tremendous boom in results and applications over the 
past 10 years 
(see references). 

As usual,  a {\em graph} is a pair $G=(V,E)$, where $V$ is the set of 
vertices (typically  $\{1,\dots,n\}$ or a subset thereof) and $E$ is 
the set of edges (an edge is a two-element subset of vertices).  
A {\em general graph} allows  multiple edges and/or loops.  
Every graph or general graph considered here  is finite (finite number
of vertices and
finite number of edges) and has a nonempty vertex set.
The {\em order} of a graph $G$, denoted $\ord{G}$, is the number of
vertices of $G$.

Let $\sRnn$ denote the set of real symmetric $n \times n$ 
matrices. For $B\in\sRnn$, the {\em graph} of $B$, denoted $\G(B)$, is the 
graph with vertices
$\{1,\dots,n \}$  and edges $\{ \{i,j \} |~b_{ij} \ne 0 \mbox{ and } i
\ne j \}$.  Note that the diagonal of $B$ is ignored in determining $\G
(B)$.   In addition, we let $\cal{S}$$(G) = \{ B \in S_{n}(\R) :$  
$\cal{G}$$(B) = G \}$. Observe that
for a given graph $G$, the classical matrices such as the {\em 
adjacency matrix} of $G$, the
{\em Laplacian matrix} of $G$, and the {\em signless Laplacian 
matrix} of $G$ 
all lie in $\cal{S}$$(G)$.

\begin{figure}[!h]
\begin{center}
\scalebox{.6}{\includegraphics{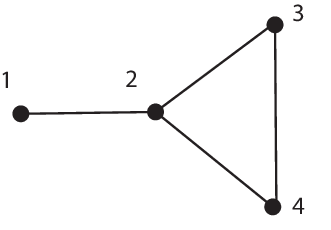}}
\caption{The graph $\G(B)$ for $B$ in
Example~\ref{ex1}}\label{fig1}
\end{center}
\end{figure}

\begin{ex} \label{ex1} {\rm 
For the matrix
$B =
\left[ \begin{array}{cccc}
0 & 1 & 0 & 0 \\
1 & 3.1 & -1.5 & 2  \\
0 & -1.5 & 1 & 1  \\
0 & 2 & 1 & 0   \end{array} \right]$, 
$\G(B)$ is shown in Figure \ref{fig1}.
}
\end{ex}

Then the {\em minimum rank} of a graph $G$ of order $n$ is defined to be
\[\mr(G)=\min\{ \rank B :~B \in \sRnn \mbox{ and }\G(B)=G \}.\]  
The problem of determining $\mr(G)$ is often referred to as the {\em standard minimum rank problem}.
The {\em maximum multiplicity} of $G$ is  given as
\[M(G) = \max\{\al_B(\lambda) :~ \lambda\in\R, B \in \sRnn 
\mbox{ and }\G(B)=G
\}.\]
Translating by a scalar matrix if necessary, it is clear that the 
maximum multiplicity of any eigenvalue is the same as maximum multiplicity of 
the eigenvalue 0. Thus maximum multiplicity is sometimes called 
{\em maximum nullity} or even {\em maximum corank}.

The following results are well-known, straightforward, and were 
presented in \cite{FH}.
\begin{enumerate}
\item 
$M(G)+\mr(G)=\ord{G}$.  

\item $\mr(G)\le \ord{G}-1$.

\item\label{Pn} $\mr(P_n)=n-1$, ($P_n$ denotes the path on $n$ vertices).

\item\label{oKn} $\mr(K_n)=1$, and if  $G$ is connected, mr($G$) = 1
implies
$G = K_{\ord G}$, that is, $G$ is the complete graph on $\ord G$ vertices.
\end{enumerate}

\begin{ex}\label{ex12}
{\rm Let  $G$ be the graph in Figure \ref{fig1} and let $A=
\left[ \begin{array}{cccc}
1 & 1 & 0 & 0 \\
1 & 2 & 1 & 1 \\
0 & 1 & 1 & 1  \\
0 & 1 & 1 & 1   \end{array} \right].$  
Since $\G(A)=G$, $G\ne K_4$, and $\rank A=2$, it follows that  $\mr(G)=2$.}
\end{ex}

For a more detailed introduction to this topic and a broad list of fundamental 
results on the minimum rank of graphs, please consult \cite{FH}. This present
survey is divided into five sections. The next section represents an update on 
recent advances and directions regarding the 
standard minimum rank problem
(that is, on the parameter $\mr(G)$). The third section of this survey 
discusses
a variant of minimum rank restricted to the subset of $\cal{S}$$(G)$ consisting 
of all the positive semidefinite matrices and the corresponding positive
semidefinite
minimum rank of a graph. The fourth section reviews a recent combinatorial
parameter, known 
as  the zero forcing  
number, and outlines  
its history and  various other
types of zero forcing 
parameters, along with their connection to the maximum nullity of a graph. 
The final two sections
are devoted to problems that are  
related to the minimum rank of a 
graph but are different in concept. 
For example, in Section \ref{sasymmr}, we consider the ranks of matrices  associated with 
directed graphs and sign patterns, and in Section \ref{sother}, we discuss related problems for (simple) graphs, such as matrices over fields other than the real numbers, the inverse inertia problem, and
minimum skew rank.

\section{Update on the standard minimum rank problem 
}\label{ssym} 

Since the first survey \cite{FH} four years ago, about 20 papers have 
appeared with results about the standard minimum rank problem, i.e., the 
problem of determining the minimum rank $\mr(G)$ of a simple graph 
describing the off-diagonal nonzero pattern of real symmetric matrices.  
A variety of results have appeared computing minimum rank for specific 
families of graphs, e.g., graphs of order at most 7 \cite{small}, 
equivalence class graphs \cite{FP}, ciclos and estrellas \cite{estrella}.   
Huang, Chang, and Yeh study various families having maximum nullity equal to zero forcing 
number 
(see Section \ref{sZ}), including block-clique graphs and unit-interval graphs \cite{Yeh}.  
Barioli, Fallat, and Smith characterized graphs having minimum rank 
equal to diameter \cite{BFS}.  Barrett et al determined the effect on 
minimum rank of certain graph operations such as edge subdivision  
\cite{BetalSubdiv}.  
Hogben and Shader studied the effect on maximum nullity of requiring null 
vectors to be generic \cite{HS}.
Recall that a path cover  of a graph $G$  is  a set of vertex disjoint induced paths that cover 
all the vertices of $G$, and the path cover number $\PC(G)$ 
is the minimum number of paths in a path cover of $G$. 
Sinkovic showed that for an outerplanar graph $G$, $\M(G)$ is bounded 
above by the path cover number $\PC(G)$  \cite{SINK}.  

Section \ref{soft} below describes computer programs that are now available 
for computing the minimum rank of small graphs.  Section \ref{savemr} 
describes work on determining the average minimum rank over all (labeled) 
graphs of a fixed order.
Some of the progress on the standard minimum rank problem is discussed 
in other sections of this article.  The zero forcing number, whose 
terminology was developed at the AIM workshop \cite{AIM06}, has played 
a role in much of the recent progress on minimum rank.  This parameter, 
its extensions, and applications to physics are described in Section 
\ref{sZ}.
The graph complement conjecture (GCC) was posed as a question at the 
AIM workshop.  Although still unproved, progress has been made on GCC, 
and it is now believed that  stronger positive semidefinite versions of 
the conjecture are true.  Thus work on GCC is discussed in Section 
\ref{spsdmr} rather than in this section.  The delta conjecture was 
also discussed at the AIM workshop  (and a stronger  
positive 
semidefinite version was conjectured by Maehara in 1987 \cite{LSS89}); 
the delta conjecture is discussed in Section \ref{spsdmr}.

\subsection{Software for minimum rank, maximum nullity, and zero 
forcing number}\label{soft} 

Since 2008 several programs have been written in the computer mathematics 
system {\em Sage} to  compute various known bounds on minimum rank and 
maximum nullity for a given graph.  For a small graph (e.g., order 
at most 10) the upper and lower bound are often equal, thereby providing 
the minimum rank.  These were originally published in  
\cite{SAGE}.  Subsequently, improvements have been made, primarily to 
the computation of zero forcing parameters, and the 2010 state of the 
art version is available at \cite{Jason}.  These programs enabled 
experimentation that led to the discovery of the estrella $S_5(K_4)$, 
a 3-connected planar graph that has the property 
$\M(S_5(K_4)^d)\ne \M(S_5(K_4))$ \cite{estrella} (here 
$G^d$ is the dual of $G$).

\subsection{Average minimum rank}\label{savemr} 
Most of the graph families for which minimum rank has been computed 
are {\em sparse}, meaning that the number of edges is much less than 
the maximum possible number of edges ${n\choose 2}\approx \frac 1 2 n^2$, 
and most are structured and exhibit symmetry. However, the random graph 
$G(n,\frac 1 2)$ for which it is equally likely that each edge is 
present or absent (the probability of each edge is $\frac 1 2$) is 
expected to have $\frac{n(n-1)} 4$ edges.  Thus the graphs for which 
minimum rank has been computed tend to present a somewhat atypical 
picture.  
Hall, Hogben, Martin, and Shader \cite{HHMS} obtained bounds on the 
average value of minimum rank  (over all labeled graphs of a fixed order). 

Formally,  the {\em average minimum rank} of graphs of order $n$ is the 
sum over all labeled graphs  of order $n$ of the minimum ranks of  the 
graphs,   divided by the number of (labeled) graphs of order $n$.  
That is,
\[\amr(n)=\frac {\sum_{\ord G = n} \mr(G)}{2^{n \choose 2}}.\]
The average  minimum rank  is equal to the expected value of the 
minimum rank of $G(n,\frac 1 2)$, denoted by $E[\mr(G(n,1/2))]$. 
The main results on average minimum rank are

\begin{thm}\label{main}
 For $n$ sufficiently large,
  \begin{enumerate}
     \item \label{it:main2} $0.146907n<\amr(n)<0.5n+\sqrt{7n\ln n}$, and
           \item \label{it:main1} 
           $\displaystyle\left|\mr(G(n,1/2))-\amr(n)\right|<\sqrt{n\ln\ln n}$ 
           with probability approaching $1$ as $n\rightarrow\infty$.
  \end{enumerate}
\end{thm}

The results in \cite{HHMS} are somewhat more general. Asymptotic bounds 
are obtained for $E[\mr(G(n,p))]$, the expected value of the minimum 
rank of $G(n,p)$, where $p$ is the probability that an edge is present, 
and for the expected value of the Colin de Verdi\`ere type parameter $\xi$.


\section{Positive semidefinite minimum rank}\label{spsdmr} 
Associating mathematical objects  
to the vertices of a graph has long been a useful tool 
in graph theory. This technique also has roots in certain minimum rank problems.

A standard example is assigning vectors to the vertices of a graph in such a
way that orthogonality corresponds to non-adjacency. That is, for any pair
of vertices $u,v$ in $G$, the vectors $\bx_u$ and $\bx_v$ assigned to $u$ and 
$v$ are orthogonal if and only if $\{u,v\} \not\in E$.

For example, if $G$ is the graph from Figure \ref{fig1}, then assigning
the standard basis vector $\bfe_1$ from $\R^2$ to vertex 1, $\bfe_2 \in \R^2$ 
to vertices 3 and 4,
and $\bfe_1+\bfe_2$ to vertex 2, is a labeling of the vertices that respects 
the condition of having nonadjacent vertices assigned to orthogonal vectors.
Also observe that if 
\[ B= \left[ \begin{array}{cccc}
\bfe_1, & \bfe_1+\bfe_2, & \bfe_2, & \bfe_2 \end{array} \right], \]
then $B$ is a $2 \times 4$ real matrix such that $B^{T}B$ is a positive 
semidefinite matrix in $\SG$.  
Moreover, the rank of $B^T B$ is two. Hence
the minimum rank among all positive semidefinite matrices in $\SG$ is at
most two (in fact, it is exactly two in this instance 
and $B^TB$ is equal to  the matrix $A$ in Example \ref{ex12}).

For any graph $G$ of order $n$, we let $\SGp$ denote the subset of $\SG$ consisting of
all real positive semidefinite matrices. Further, we let
\[ \mrp(G) = \min \{ \rank A : A \in \SGp \}, \]
and
\[ \Mp(G) = \max \{  \nul A : A \in \SGp \}. \] The parameter $\mrp(G)$ is
called the {\em (real) minimum positive semidefinite rank of $G$}, while
$\Mp(G)$ is called the {\em maximum positive semidefinite nullity of $G$}.
As with the case of standard minimum rank, it is clear that for any graph $G$
\[ \mrp(G) + \Mp(G) = \ord{G}.\]

Now, following the example above, if $G$ is a graph and for each vertex 
$i \in V$ we assign the vector $\bv_i \in \R^d$ such that 
$\bv_i^{T} \bv_j =0$ if and 
only if $\{i,j\} \not\in E$, then the matrix
$B^T B$, where $B=[\bv_1, \bv_2, \ldots, \bv_n]$ is in $\SGp$ with rank equal to $k$.
Such a vector representation is called an
{\rm orthogonal vector representation} (see also \cite{LSS89, LSS00} where 
the representation above is known as a {\em faithful} 
orthogonal vector representation). Orthogonal vector representations also 
arise in the works \cite{Netal08chord, Netal09OS, MNZ} mostly over the 
complex field, but the concept is analogous. Orthogonal vector representations
(of the non-faithful variety) also appear in connection with the 
Lov\'asz $\vartheta$ function and related versions of certain sandwich 
type theorems (see, for example, \cite{FH} and the relevant references within). 
It follows easily that $\mrp(G)$ coincides with the smallest $d$ such that 
$G$ admits an orthogonal vector representation with vectors from $\R^d$.

As noted above, it is also of interest to investigate the smallest $d$ such
that the graph $G$ admits an orthogonal vector representation with vectors lying 
in $\C^d$ instead of restricting to the real case. The smallest such $d$ will
be denoted by $\mrp^{\C}(G)$ and it is not difficult to observe that
\[ \mrp^{\C}(G) = \min \{ \rank A : A \in \PSDG\},\]
where $\PSDG$ is the subset of positive semidefinite matrices among all
complex Hermitian matrices $A$ such that $\G(A)=G$. This term has been 
well-studied just like its real counterpart and in the papers 
\cite{Netal08chord, Netal09OS, MNZ} we note that $\mrp^{\C}(G)$ is denoted
by the symbol msr$(G)$. It is very important to observe that changing 
fields from $\R$ to $\C$ does result in a different parameter as noted in
\cite{smallparam}. 

In many ways, it does appear that the parameters $\mrp(G)$ and
$\mrp^{\C}(G)$ may be more natural graph-type parameters when compared to 
other notions of minimum rank. This opinion may be defended by the 
simplicity of many results about minimum positive semidefinite rank and
its connections to graph theory. 

For example, it is known that the minimum semidefinite rank of any tree is
precisely the order of the tree less one, which is as large as the minimum
semidefinite rank can be in general (see, for example, \cite{vdH09} or
\cite{Netal08chord}).

In the context of certain graph operations, the minimum semidefinite rank
behaves rather nicely. For example, in the case when $G$ has a cut vertex
the minimum semidefinite rank of $G$ can be computed by summing the 
minimum semidefinite ranks of smaller graphs (see \cite{Netal08chord} for a proof over 
the complex numbers, although a similar argument will work over the reals, see also \cite{vdH09}). We
note here that the formula below can easily be used with a simple induction
argument to verify that the minimum semidefinite rank of trees is precisely
the order of the tree less one.

\begin{thm}\label{cut} 
Suppose $G$ has a cut-vertex $v$. For $i=1,\dots, h$, 
let $W_i\subseteq V(G)$ be the vertices of the $i$th component 
of $G-v$ and let $G_i$ be the subgraph induced by $\{v\}\cup W_i$. Then
 \[
   \mrp(G)=\sum_1^h \mrp(G_i).
\]
\end{thm}

 An analogous cut-vertex reduction formula for 
minimum rank for a graph was obtained earlier by various authors and is presented in \cite{FH}. However
that formula is more convoluted and depends on the notion of the rank-spread of
a vertex. Recall that the
rank-spread of $G$ at vertex $v$ is defined to be
$r_{v}(G)=\mr(G)-\mr(G-v).$ In the positive semidefinite case it is not
difficult to observe that the rank spread of a vertex $v$ is bounded between
\[ 0 \leq \mrp(G)-\mrp(G-v) \leq {\rm deg}(v), \]
where ${\rm deg}(v)$ is the degree of the vertex $v$. The fact that the rank
spread in the positive semidefinite case can be larger than 2
seems to simplify
calculations in the case of cut vertex reduction.

In the case of the join of two graphs a similar simplification occurs. 
Recall that the  {\em join} $G\vee G'$ of two disjoint graphs $G=(V,E)$ and
$G'=(V',E')$ is the union of $G\cup G'$ and the complete
bipartite graph with with vertex set $V\cup V'$ and partition
$\{V,V'\}$. 

The following fact was proved in \cite{Netal09OS} over the complex numbers and 
in \cite{GCC} over the reals.

\begin{thm} \label{join}
If $G$ and $H$ are two graphs, then 
\begin{equation}
\mrp(G \vee H) = \max \{ \mrp(G \vee K_1), \mrp(H \vee K_1)\},\label{mrp-join}
\end{equation}
where $K_1$ is the complete graph on a single vertex.
\end{thm}

Observe that if $G$ and $H$ do not contain any isolated vertices, then 
we have \[ \mrp(G \vee H) = \max \{ \mrp(G), \mrp(H)\}.\]

For standard minimum rank it is well-known that the equations above need
not hold in general, and, in fact, $\mr(G \vee H)$ only behaves as in 
(\ref{mrp-join}) for the special
case of graphs that are among the so-called inertia-balanced and not anomalous (see \cite{FH}). For example,
if $G$ and $H$ are both trees or are both decomposable graphs, then (\ref{mrp-join}) is valid.

In addition, many other facts are known about the parameter $\mrp(G)$. For 
example, it is easy to verify that for any graph $G$, $\mr(G) \leq \mrp(G) \leq 
\cc(G)$, where $\cc(G)$ denotes the clique cover number of $G$ (that is, the 
fewest number of cliques needed to cover the edges of $G$). Furthermore, if $G$
is known to be chordal (no induced cycles of length four or more), then 
$\mrp(G) = \cc(G)$. (See \cite{Netal08chord} for a proof over the complex numbers. This equation over
the reals then follows easily.) However, 
$\mr(G) < \cc(G)$ for 
 any chordal graph for which it is known that $\mr(G)<\mrp(G)$, such as a tree that is not a path.

In addition, many other interesting facts are known about the minimum semidefinite rank,
including:
\begin{itemize}
\item If $G'$ is obtained from $G$ by an edge subdivision, then
$\mrp(G')=\mrp(G)+1$ (see \cite{JMN}, a similar argument applies over $\R$),
\item  If $G$ is triangle free, then $\mrp(G) \geq \mrp^{\C}(G) \geq  \lceil n/2 \rceil$ (see 
\cite{Dthesis, Dea}),
\item If $G$ is outerplanar, then $\Mp(G)$ is equal to the tree cover number of 
$G$ (see \cite{BFMN}).\end{itemize}

\subsection{Delta conjecture} 
As mentioned in Section 2,  
at the AIM workshop in 2006 an interesting
inequality was conjectured to hold between the minimum degree and maximum 
nullity (see \cite{AIMrpt}). Since that time 
the validity of this inequality is still unresolved.
However, there is significant positive evidence to suggest that the inequality
is indeed valid. The {\em delta conjecture}, as it has become known, states
that any graph $G$ with minimum degree $\delta(G)$ satisfies,
\[ \M(G) \geq \delta(G).\] Equivalently, we 
could ask if $\mr(G) \leq \ord{G} -\delta(G)$ holds for all graphs $G$. 

At present, the delta conjecture is known to hold for many classes of graphs
including trees, graphs with $\delta(G) \leq 3$, 
bipartite graphs (see \cite{delta}), along with
various other examples. 

A stronger version of the delta conjecture involving 
positive semidefinite matrices has also been suggested and at present
remains open (see also \cite{LSS89} 
 for a reference to a conjecture made
by Maehara). Is it true that for all graphs $G$, $\Mp(G) \geq \delta(G)$?
If this inequality holds, then the delta conjecture would be solved, as
$\M(G) \geq \Mp(G)$. However, at present the relationship between $\Mp(G)$ and 
$\delta(G)$ has not been fleshed 
out, and still remains for the most part
a mystery. For example, it is not known if  $\Mp(G) \geq \delta(G)$ for
bipartite graphs $G$.

On the other hand, there
is a nice connection between $\Mp(G)$ and the vertex connectivity of a graph, 
denoted by $\kappa(G)$.  In \cite{LSS89}, it was shown that $\Mp(G)\ge\kappa(G)$.  Unfortunately, it is also known that $\delta(G)\ge\kappa(G)$ and strict inequality is possible.
Recall that the Colin de Verdi\`{e}re parameter $\nu(G)$ (see \cite{FH} for a basic introduction on
this topic)  
is defined to be the maximum 
multiplicity of 0 as an eigenvalue
among matrices $A\in\sRnn$ that satisfy:
\begin{itemize}
\item $\G(A )=G$.
\item $A$ is positive semidefinite.
\item $A$ satisfies the Strong Arnold Hypothesis.
\end{itemize}
In \cite{vdH08} is was observed that results in \cite{LSS89} in fact implied that $\nu(G)\ge\kappa(G)$.  In \cite{param} it is conjectured that $\nu(G)\ge\delta(G)$.

\subsection{Graph Complement Conjecture (GCC)} 
Another interesting conjecture that arose from the 2006 AIM workshop has
become known as the graph complement conjecture or GCC for short
(see \cite{AIMrpt}). The GCC can be written as the following
conjecture  about
the minimum rank of $G$ and its complement,
\begin{equation}
\mr(G) + \mr(\Gc) \leq \ord{G}+2,\label{gcc-eq}
\end{equation} where $\Gc$ is the complement of $G$.

For instance, if $G=C_5$, the cycle on 5 vertices, 
then $\mr(C_5)=3$ and $\mr(\overline{C_5})=\mr(C_5)=3$. Hence, 
$\mr(G) + \mr(\Gc)= 3+3 < 5+2$. For paths
on $n$ vertices, it can be shown that equality holds in (\ref{gcc-eq})
whenever $n \geq 4$ (see \cite{AIM}). 

As with the delta conjecture, there is overwhelming evidence in favor of GCC,
however it remains unresolved at present. In addition, stronger forms of
GCC have since been suspected and remain open. For example, is
the inequality \[ \mrp(G) + \mrp(\Gc) \leq \ord{G}+2 \] valid in general?

Observe that GCC (and it variants) can also be stated equivalently in terms
of maximum nullities. For example,
\[ \M(G) + \M(\Gc) \geq \ord{G}-2, \; {\rm and} \; 
\M_{+}(G) + \M_{+}(\Gc) \geq \ord{G}-2.\] 
A further strengthening has also been conjectured in terms of the Colin de 
Verdi\`{e}re 
parameter $\nu(G)$ (see \cite{GCC}):
\[ \nu(G) + \nu(\Gc) \geq \ord{G}-2.\]

In the recent work \cite{GCC} there is a number of positive results pertaining
to  the GCC and it variants, including the case of the join of two 
graphs and restrictions to $k$-trees.


\section{Zero forcing parameters}\label{sZ} 

One approach to studying the minimum rank or maximum nullity of a graph is to
investigate the possible structure of the null space in 
order to provide bounds on the nullity itself. 

For example, if the null space of a given $n \times n$ 
matrix $A$ has dimension at least 
2 (or $> 1$), then for each $i=1,2,\ldots, n$, there exists a nonzero 
vector $\bx$ in 
the null space of $A$ with $x_i =0$. Another way to view this concept is the 
following: Suppose there exists an index $i$ such that any null vector $\bx$ with
$x_i =0$ implies $\bx=0$. Then we may conclude that the dimension of $\nul A$ 
cannot be more than one. More generally, if a set $S$ of indices has the 
property that  $A\bx=0$ and $x_i=0$ for all $i\in S$ implies $\bx=0$, 
then $\nul A\le \ord S$.

Consider the path on $n$ vertices as a preliminary example. Suppose 
$A \in \mathcal{S}(P_n)$, and that the vertices of $P_n$ are labeled
in increasing order. Suppose  that $\bx$ is a null vector for $A$ and that
$x_1=0$. Then the equation $A\bx=0$ in the first coordinate becomes
\[ a_{11}x_1 + \sum_{1 \sim j} a_{1j} x_j = a_{11}x_1 + a_{12}x_2 =0,\]
where $i \sim j$ means vertex $i$ is adjacent to vertex $j$. The above 
equations imply that $x_2=0$ as $a_{12} \neq 0$. Replacing $i=1$ with $i=2$ and
continuing in the same manner we deduce that $x_3=0$. In other words, if 
$A \in \mathcal{S}(P_n)$, then the dimension of $\nul A$ is at most 1. Hence
we may conclude that $\M(P_n)=1$. 

More generally, if $A \in \SG$, then for each $i$ the $i$th coordinate of the 
equation $A\bx=0$ may be written as
\begin{equation}
 a_{ii}x_i + \sum_{i \sim j} a_{ij} x_j =0.\label{zfs-eq} \end{equation}
Appeal to (\ref{zfs-eq}) to provide some intuition as to when a collection of
zero coordinates in a null vector of $A$ necessarily implies that the null
vector must have been the zero vector to start with. For instance, suppose
$x_i=0$ and $x_j=0$ for all but one neighbor of $i$. Then by (\ref{zfs-eq}),
we have that all of the neighbors of $i$ will have zero coordinates
in $\bx$. If this process could continue to demonstrate that $\bx=0$, then we may
conclude that the dimension of $\nul A$ cannot exceed the number of 
neighbors of $i$. To formalize this idea, we devise a coloring scheme on the 
vertices of $G$.

Suppose $G = (V,E)$ is a given graph and that the vertices of $G$  
are partitioned into two sets, $V = B \cup W$, where the vertices in $B$ are
colored black and the vertices in $W$ are colored white. The 
goal of the game is
to color all of the vertices in $G$ black. To do this, we define a rule
known as a {\em color change rule}. The color change rule in 
this case, denoted by CCR-Z, is as follows: a black vertex $v$ 
can color a white
neighbor $u$ if it is the only such white neighbor of $v$. In this case, we 
say that $v$ {\em forces} $u$. The rule corresponds to the implication that we 
observed above in (\ref{zfs-eq}), if we associate the black vertices in $B$
with the initial zero coordinates of a given null vector. 

\begin{figure}[!h]
\begin{center}
\scalebox{.6}{\includegraphics{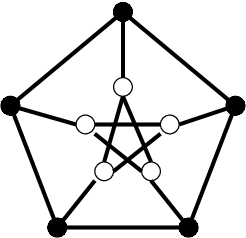}}
\caption{Petersen graph}\label{petersen}
\end{center}
\end{figure}

Furthermore, a subset of vertices $S \subset V$ is called a {\em zero forcing 
set for $G$} if whenever the vertices of $S$ are colored black while all
remaining all colored white, then all vertices of $V$ are forced to be black
under repeated application of the color change rule CCR-Z. 
For example, a pendent vertex of a path is
a zero forcing set for that path. If $G$ is the Petersen graph shown in Figure
\ref{petersen}, then the vertices colored black form a zero forcing set.

In other words a zero forcing set of 
vertices corresponds to an initial collection of indices with the property
that if the coordinates of these indices are assigned with zeros in a
null vector, then the 
associated null vector must be the zero vector.

A subset of the vertices is called a {\em minimum zero forcing set for G} if it
is a zero forcing set for $G$ and there is no other zero forcing sets that
consist of fewer vertices. For example, a pendant vertex of a path is a 
minimum zero forcing set of a path, and the set of five black vertices in
the Petersen graph above form a minimum zero forcing set for the Petersen
graph. Finally, the size of a minimum zero forcing set for $G$ is called 
the {\em zero forcing number of $G$}, and is denoted by $\ZFS(G)$ \cite{AIM}. Thus
$\ZFS(P_n)=1$ and the zero forcing number of the Petersen graph  
is 5.

From the construction of $\ZFS(G)$ it follows that the zero forcing number of 
a graph is always an upper bound for the maximum nullity of $G$.

\begin{thm} 
{\rm \cite{AIM}}
For any graph $G$, $\ZFS(G) \geq \M(G)$.
\label{zfs-thm}
\end{thm}

Other  
properties of $\ZFS(G)$ can be found in \cite{smallparam} including
non-uniqueness of minimum zero forcing sets, the intersection over all 
minimum zero forcing sets of a graph is always empty, 
and $\ZFS(G)$ is 
always an upper bound for $\PC(G)$ (the path cover number of $G$). 
Other properties of $\ZFS(G)$ can be found in the works \cite{EHHLR, Owe, Row}.

The case of equality between $\ZFS(G)$ and $\M(G)$ is still very much of
interest and unresolved, namely 
the problem of 
characterizing the graphs $G$ for 
which $\M(G) = \ZFS(G)$. For example, equality holds between these two
parameters for trees and various other  
examples \cite{AIM, Yeh}. However, it is known
that the gap between $\ZFS(G)$ and $\M(G)$ can grow without bound on a 
sequence of graphs.

 The idea of zero forcing on a graph was introduced independently by physicists to study control of quantum systems \cite{graphinfect, Sev}.  Vertices are colored black or white and the same color change rule is applied but the process is   called  {\em propagation}, and again it is the minimum number of vertices in a set that propagates (i.e., the zero forcing number $\ZFS(G)$) that is of interest.   Quoting from \cite{graphinfect}, ``Our goal is to
determine if such a configuration [of black vertices] is compatible with being
a nontrivial eigenstate of the network Hamiltonian (that is
an eigenstate in which not all the vertices are black); if not,
then the whole network can be controlled."  in this context, ``nontrivial" refers to non-constant spin orientation.  There seem to be deep connections between minimum rank problems and control of quantum systems that are only beginning to be explored. 

To bound the maximum nullity of  different sets of matrices described by a graph (such as positive definite matrices), variations of zero forcing 
have been defined by varying the color change rule as needed. Given a color change rule CCR-$x$ and a coloring of of a graph $G$, the {\em derived set} is the 
set of black vertices obtained by applying CCR-$x$ 
until no more changes are possible.
A {\em (CCR-$x$) zero forcing set} for  $G$ is a subset of vertices $Z$ 
such that if initially the vertices in $Z$ are colored black 
and the remaining vertices are colored white, then the derived 
set  is  all  the vertices of $G$. 
The {\em (CCR-$x$) zero forcing number}   is the minimum of $\ord Z$ 
over all (CCR-$x$) zero forcing sets $Z\subseteq V(G)$. 

\subsection{Positive semidefinite zero forcing}\label{sspsdZ}
The analogous concept of zero forcing in the positive semidefinite case comes with its own version of a color change rule.
The {\em positive semidefinite color change rule} \cite{smallparam} 
 is:
\bit
\item [CCR-$\Zp$] Let $B$ be the set consisting of all the black vertices of $G$.  
Let $W_1,\dots, W_k$ be the sets of vertices of the $k$ components of
 $G-B$ (note that it is possible that $k=1$). 
  Let $w\in W_i$.
If $u\in B$ and $w$ is the only white neighbor of  $u$ in $G[W_i\cup B]$, 
then change the color of $w$ to black. 
\eit As indicated above, the
{\em  positive semidefinite zero forcing number of a  graph $G$}, 
denoted by $\Zp(G)$,   is the minimum of $\ord X$ over all CCR-$\Zp$ zero forcing sets $X\subseteq V_G$.  

Forcing with the positive semidefinite color change rule
can be viewed as decomposing the graph into a union of certain  
induced subgraphs and then using CCR-Z on each of these induced subgraphs.
For example, it is evident that $\Zp(T)=1$ for any tree $T$, because any one
vertex is a positive semidefinite zero forcing set for $T$. In addition, it
is also easy to verify that $\Zp(G) \leq \ZFS(G)$ for any graph $G$.

\begin{figure}[!h]
\begin{center}
\scalebox{.6}{\includegraphics{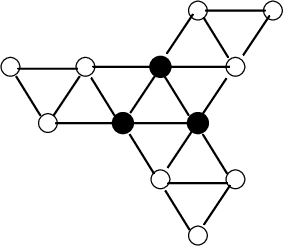}}
\caption{The Pinwheel on 12 vertices}\label{pinwheel}
\end{center}
\end{figure}

The graph $G$ in Figure \ref{pinwheel}
satisfies $\Zp(G) =3 < 4 =\ZFS(G)$
\cite{smallparam}; 
the vertices colored black form a minimum positive 
semidefinite zero forcing set.

As with the case of standard zero forcing, the positive semidefinite
zero forcing number is always an upper bound on the positive 
semidefinite maximum nullity.

\begin{thm} 
{\rm \cite{smallparam}}
For any graph $G$, $\Zp(G) \geq \Mp(G)$.
\label{zfspsd-thm}
\end{thm}

We also note 
 that the concept of 
 positive semidefinite zero forcing is related to the notion of
ordered sets that appear in \cite{Netal09OS, JMN, MNZ}. In fact, it is known 
(see \cite{smallparam}) that for any graph $G=(V,E)$ and any ordered set $S$, 
$V\setminus S$ is a positive semidefinite forcing set for $G$, 
and for any positive semidefinite forcing set $X$ for $G$, 
there is an order that makes $V\setminus X$ an ordered set for $G$.  
Thus $\Zp(G) + OS(G)=|G|$ (here $OS(G)$ is the ordered set number of $G$,
see \cite{Netal09OS}). It is also known that $\ZFS(G)$ is related in a similar
manner to the connected ordered set number \cite{MNZ}.

From the relation $\Zp(G) + OS(G)=|G|$ and the 
fact that $OS(G) \leq |G| -\delta(G)$ from \cite{MNZ},  for any graph $G$ we have 
\[ \Zp(G) \geq \delta(G), \]

\subsection{Other zero forcing parameters}\label{ssoZ}

In an attempt to obtain improved bounds on $\M$,  
graphs that allow loops have been considered \cite{param}.
A {\em loop graph} is a graph that allows loops, i.e., 
$\wh G =(V_{\wh G},E_{\wh G})$
where $V_{\wh G}$ is the set of vertices of $\wh G$ and the set of 
edges $E_{\wh G}$
is a set of two-element multisets.  Vertex $u$ is a {\em neighbor}
of vertex $v$
in $\wh G$ if $\{u,v\}\in E_{\wh G}$;
note that $u$ is a neighbor of itself if and only if the loop 
$\{u,u\}$ is an edge.
The {\em underlying simple graph} of a loop graph $\wh G$ is 
the graph $G$ obtained
from $\wh G$ by deleting all loops.  

The {\em set of symmetric matrices described by a loop graph $\wh G$} is
\[\SGh=\{A=[a_{ij}] \in\sRnn : a_{ij}\ne 0
\mbox{ if and only if } \{i,j\}\in E_{\wh G}\}.\]
Note that a loop graph $\wh G$ constrains the zero-nonzero pattern 
of the main diagonal
entries of matrices described by $\wh G$.
There is a distinction between a  graph, i.e., a simple graph,
and a loop graph that has no loops---the latter forces the matrices 
to have zero diagonal,
whereas the former does not 
(see also Section \ref{ssasymparam}).
The color change rule for loop graphs is:
\bit
\item [CCR-$\ZFS(\wh G)$] If exactly one  neighbor $w$ of $u$  is white, then change the color of $w$ to black.\eit 
 The  {\em   zero forcing number of a loop graph $\wh G$}, denoted by $\ZFS(\wh G)$,   is the zero forcing parameter for CCR-$\ZFS(\wh G)$.
The {\em  enhanced zero forcing number} 
of a  
(simple) graph $G$, denoted by $\Zh(G)$, is the maximum of $\ZFS(\wh G)$ 
over all   loop graphs $\wh G$ such that the underlying simple graph of 
$\wh G$ is $G$ (see \cite{param}).

\begin{thm} 
{\rm \cite{param}}
\label{propZhat} For any graph $G$, $\M(G) \leq\Zh(G)\le \ZFS(G)$.
\end{thm}

Finally, the {\em  loop zero forcing number} 
of a 
(simple) graph $G$, denoted by $\Zl(G)$, is  $\ZFS(\wh G)$ where $\wh G$ 
is the 
specific 
loop graph whose underlying simple graph is $G$, and such that 
$\wh G$ has a loop at $v\in V_G$ if and only if $\deg_G v\ge 1$.

Although $\Zl$ is already defined as $\ZFS$ evaluated on a specific 
loop graph, we can see that $\Zl$ is a zero forcing parameter, which 
aids in computing the value of this parameter.  
The color change rule associated with the loop zero forcing number 
is  (see also \cite{param}): 
\bit
\item[CCR-$\Zl$] If $u$ is black and exactly one  neighbor $w$ of $u$  is white, 
then change the color of $w$ to black.  If $w$ is white, $w$ has a 
neighbor, and every  neighbor of $w$ is black, then change the color 
of $w$ to black.
\eit

\begin{thm} 
{\rm \cite{param}} \label{propZhat2} For any graph $G$, 
$\Zp(G)\le \Zl(G) \leq\Zh(G)$. 
\end{thm}

\vspace{-0.5cm}
\begin{figure}[!h]
\begin{center}
\caption{Relationships between zero forcing parameters, parameters 
related to maximum nullity, and other graph parameters.}\label{pardiag}
\vspace{1.1 cm}
\scalebox{.85}{
\begin{picture}(140,300)(-25,-100)

\put(50,236){\small $\ZFS$} 
\put(53,234){\line(0,-1){33}} 
\put(53,234){\line(1,-6){45}} 
\put(94,-45){\small  $\PC$}

\put(50,191){\small  $\Zh$}
\put(50,190){\line(-2,-3){21}} 
\put(53,190){\line(0,-1){55}} 

\put(22,150){\small  $\Zl$}
\put(22,148){\line(-2,-3){15}} 

\put(-1,121){\small $\Zp$} 
\put(3,119){\line(0,-1){48}} 
\put(2,119){\line(-2,-3){25}} 

\put(-32,75){\small $\tw$}
\put(-27,73){\line(0,-1){109}} 

\put(49,129){\small $\M$}
\put(48.5,126){\line(-3,-4){40}} 
\put(53,126){\line(0,-1){47}} 

\put(-1,63){\small $\Mp$}
\put(3,60){\line(0,-1){48}} 

\put(50,71){$\xi$}
\put(50,68){\line(-3,-4){43}} 
\put(53,68){\line(0,-1){55}} 
\put(50,7){\small $\mu$}

\put(1,5){$\nu$}
\put(3,3){\line(-1,-3){30}} 
\put(6,5){\line(-4,-6){30}} 

\put(-30,-43){$\delta$}
\put(-28,-45){\line(0,-1){40}} 
\put(-30,-93){$\kappa$}

\end{picture}
}
\end{center}
\end{figure}

Figure \ref{pardiag}, 
adapted from \cite{param},
describes the relationships between the zero forcing parameters, 
maximum nullity parameters, and other graph parameters (for graphs 
that have at least one edge). In Figure \ref{pardiag}, a line between 
two parameters $q,p$ means that for all graphs 
$G$, $q(G)\le p(G)$, where $q$ is below $p$ in 
the diagram. Furthermore, it is known in all cases that inequalities 
represented in Figure \ref{pardiag} can be strict (see \cite{param}).  
The strongest form of the delta conjecture ($\delta(G)\le\nu(G)$) 
appears as a dashed line of small triangles. (The parameters 
$\mu$ and $\xi$ are Colin de Verdi\`{e}re type parameters and 
are defined in \cite{FH}; $\tw(G)$ denotes the tree-width of $G$ 
(see \cite{param}).)


\section{Minimum rank of patterns and other types of graphs}\label{sasymmr}

The families of matrices discussed in previous sections have had 
off-diagonal nonzero patterns  described by edges of simple undirected 
graphs.   In this section we survey work on the minimum rank of matrices 
with more general patterns of nonzero entries, sometimes eliminating the 
requirement of positional symmetry by using directed graphs,  allowing 
the pattern to (more fully) constrain the diagonal, and including sign 
patterns in addition to nonzero patterns.  
  The minimum rank problem for  nonzero patterns has been studied over 
  fields other than the real numbers, but for simplicity we limit the 
  discussion here to matrices over the real numbers.

A {\em nonzero pattern} 
is an $m\x n$ matrix $Y$  whose entries are
elements of
$\{*, 0\}$. 
For $B=[b_{ij}]\in \R^{m\x n}$,  the {\em pattern of $B$}, 
$\Y(B)=[y_{ij}]$, is the $m\x n$ nonzero pattern with $y_{ij} = *$ if 
$ b_{ij}\ne 0$ and $y_{ij} =0$ if $b_{ij}= 0$.
 A {\em sign pattern}  is a  matrix having entries in $\{+,-,0\}$. 
 For  $B\in\R^{m\x n}$, $\sgn(B)$ is the sign pattern having entries 
 that are the signs of the corresponding entries in $B$.  An $n\x n$ 
 (nonzero or sign) pattern  is called {\em square}. 

 The definitions of 
{\em minimum rank} and {\em maximum nullity} are also extended to an 
$m\x n$ nonzero pattern or sign pattern.  
For a  nonzero pattern $Y$:  
\bea\mr(Y)&=&\min\{ \rank (B) : \, B\in\R^{m\x n}, \ \Y(B)=Y\}.\\
\M(Y) &=& \max\{\nul(B): \, B\in\R^{m\x n}, \ \Y(B)=Y\}.\eea
For a sign pattern $Y$, replace $\Y(B)=Y$ by $\sgn(Y)=Y$.
If $Y$ is $m\x n$, then
$\mr(Y)+\M(Y)=n.$  
The problem of determining the minimum rank  of a sign pattern, also 
called the {\em sign-rank}, has important applications to communication 
complexity 
(see Section \ref{ssSP}).

\subsection{Parameters related to minimum rank of nonzero patterns}
\label{ssznz} 

In \cite{square1} it is shown that the minimum rank problem for a 
nonzero pattern can be converted to a (larger) minimum rank problem 
of standard type, i.e., symmetric matrices described by a simple 
undirected graph.

A {\em $t$-triangle} of an $m\x n$ nonzero pattern $Y$ is a $t\x t$ 
subpattern that is permutation similar to a pattern that is upper 
triangular with all diagonal entries nonzero.  
The {\em triangle number} of pattern $Y$, denoted $\tri(Y)$,  is 
the maximum size of a triangle in $Y$.  The triangle  number and 
$t$-triangles have been used as a lower bound for minimum rank in 
both the symmetric and asymmetric minimum rank problems, see e.g., 
\cite{BHL04}, \cite{CJ}.
The triangle  number was a focus of the papers \cite{CJ, JL08}, 
where it was denoted MT$(Y)$.  Small patterns $Y$ for which 
$\mr(Y)=\tri(Y)$ were determined; this includes all $m\x n$ patterns 
with $m\le 5$ (the smallest known example where 
$\mr(Y)>\tri(Y)$ is $7 \x 7$).

For a square nonzero pattern $Y$, the {\em (row) edit distance to nonsingularity}, 
$\ED(Y)$, of $Y$ is the minimum number of rows that must be changed to 
obtain a pattern that requires nonsingularity {\rm \cite{square1}}.  
The edit distance to nonsingularity is related to the triangle number.

\begin{thm}\label{triED}{\rm \cite{square1}} For an $n\x n$ nonzero 
pattern $Y$, $\tri(Y)+ \ED(Y)=n$.
\end{thm}

\subsection{Graphs of various types}\label{ssasymparam}

 Graphs continue to be a powerful tool in the study of minimum rank 
 of nonzero patterns, 
 but the expansion of the type of pattern discussed necessitates 
 being more inclusive  in our definition of ``graph."  Throughout 
the remainder of
Section \ref{sasymmr}, a graph can be simple or 
 allow loops, and can be undirected or directed.
When describing a specific type of graph, we always use one of 
the terms {\em simple} or {\em loop} and one of the terms  
{\em graph} or {\em digraph}.  
We use the  term {\em graph of any type} to mean one of a simple graph, 
a loop graph, a simple digraph, or a loop digraph.  
We continue to require symmetric matrices for an (undirected)  
graph (simple or having loops), so in case this restriction 
is not desired, a doubly directed digraph (simple or having loops) 
should be used if the pattern of nonzero 
entries is symmetric. 
 Note that loop graphs were already introduced in Section \ref{ssoZ} (where a loop graph was denoted by $\wh G$), and the definitions given in that section for the set of matrices described by the graph, minimum rank, maximum nullity, zero forcing number, etc. coincide with those given here, although the notation is slightly different.

Each type of graph describes a set of matrices,  the 
{\em qualitative class} of $G$ of order $n$, denoted by $\Q(G)$.  
\bit \item For a simple graph $G$, 
$\Q(G)=\{A\in \Rnn:A^T=A \mbox{ and for } i\ne j, \, a_{ij}\ne 0 
\Leftrightarrow \{i,j\}\in E(G)\}.$

\vspace{-2mm}
\item For a simple digraph $G$, 
$\Q(G)=\{A\in \Rnn: \mbox{for } i\ne j,  \, a_{ij}\ne 0 
\Leftrightarrow (i,j)\in E(G)\}.$

\vspace{-2mm}
\item For a loop graph $G$,
$\Q(G)=\{A\in \Rnn:A^T=A \mbox{ and } a_{ij}\ne 0 
\Leftrightarrow \{i,j\}\in E(G)\}.$

\vspace{-2mm} 
\item For a loop digraph $G$, 
$\Q(G)=\{A\in \Rnn:a_{ij}\ne 0 \Leftrightarrow (i,j)\in E(G)\}.$
\eit
For a graph $G$ of any type,
\[\mr(G)=\min\{\rank A : A\in\Q(G)\}\qquad\mbox{and}\qquad 
\M(G)=\max\{\nul A : A\in\Q(G)\}.\]
Clearly $\mr(G)+\M(G)=\ord G.$

The definition of zero forcing number has been be extended from 
simple graphs to loop graphs, loop digraphs, and simple digraphs 
\cite{square1, cancun}.  In this section, we denote a graph by $G$ even if it is a loop graph (or digraph), and the zero forcing number of $G$ is denoted $\ZFS(G)$.  As noted in Section \ref{ssoZ}, the only change needed in the definition of zero forcing number is the color change rule, 
which depends on the type of graph.  The color change rules for a simple graph and a loop graph are CCR-$\ZFS$  and  CCR-$\ZFS(\wh G)$, respectively, defined in Section \ref{sZ}.  For simple and loop digraphs, the color change rules  are:
\bit 
\vspace{-2mm}\item[CCR-$\ZFS(\D)$]
Let $G$ be a   a simple digraph. If $u$ is a black vertex and exactly one  
out-neighbor $v$ of $u$ is white, then change the color of $v$ to black. 
\vspace{-2mm}\item[CCR-$\ZFS(\wh\D)$]
Let $G$ be a   a loop digraph. If  exactly one  out-neighbor $v$ of $u$ 
is white, then change the color of $v$ to black (the possibility 
that $u=v$ is permitted). 
\eit
\vspace{-2mm}
 Examples of zero forcing on various types of graphs are given in \cite{cancun}.
 Regardless of the type of graph, the zero 
forcing number bounds maximum nullity from above.

\begin{thm}\label{thmMZ}{\rm \cite{cancun}} If $G$ is any type of graph, 
then $\M(G) \le \ZFS(G)$.
\end{thm}

If $G$ is a loop digraph, the nonzero pattern of $G$ is $\Y(G)=\Y(B)$ 
where $B\in\Q(G)$, the triangle number of $G$ is  $\tri(G)=\tri(\Y(G))$, 
and the edit distance of $G$ is $\ED(G)=\ED(\Y(G))$.  These parameters  
are related.

\begin{thm}\label{thmZtri}{\rm \cite{square1}} If $G$ is a loop digraph, 
then $\tri(G)+ \ZFS(G)=\ord G$ and $\ED(G)=\ZFS(G)$.
\end{thm}

\subsection{Minimum rank of sign patterns}\label{ssSP} 

The minimum rank of 
full 
sign patterns has important applications to 
communication complexity in computer science (a sign pattern is 
{\em full} if all entries are nonzero), and  
significant 
 progress on minimum 
rank of 
full 
sign patterns has been obtained 
through work on communication complexity.  

In a simple model of communication,  
described in  \cite{CdCM}, 
there are two processors $A$ and $B$, each of which receives its 
own input (a string of bits that are $0$ or $1$), and  the goal is 
to compute a value that is a function of both inputs. The 
computation function  can be described by a $\{0, 1\}$-matrix 
$M$ with rows indexed by the possible inputs of $A$,  
columns 
indexed by the possible inputs for $B$, and the entry representing 
the value computed.  A {\em (deterministic) protocol} tells the 
processors how to exchange information to enable this computation,  
The {\em (deterministic) communication complexity} $c(M)$ associated 
to the $\{0, 1\}$ function matrix $M$ is the minimum number of bits 
that must be transmitted in any protocol associated with $M$.  Melhorn 
and Schmidt \cite{MS82} showed that 
$ \log_2\rank M\le c(M)\le \rank M$ \cite{CdCM}.  

Communication complexity is also studied from a probabilistic point 
of view
; this approach is described in \cite{Lokam}.  An {\em unbounded error probabilistic protocol} tells the 
processors how to exchange information to enable computation that will 
be accurate with probability $>\frac 1 2$.  The {\em unbounded error 
probabilistic  communication complexity} $\uppcc(M)$ associated to 
the function matrix $M$ is the minimum number of bits that must be 
transmitted in any unbounded error probabilistic protocol associated 
with $M$. When studying $\uppcc$, it is common to use a  
$\{+1,-1\}$-matrix.
A $\{0,1\}$-matrix $M$ can be converted to a $\{+1,-1\}$-matrix by 
replacing entry $m_{ij}$ by $(-1)^{m_{ij}}$, or equivalently, 
using $J-2M$, where $J$ is the all ones matrix.  If $M$ is an 
$m\x n$ $\{+1,-1\}$-matrix, then $\sgn(M)$ is a full sign pattern, 
and if $X$ is an $m\x n$ $\{+,-\}$ sign pattern, then $M_X$ denotes 
the $m\x n$  $\{+1,-1\}$-matrix having $\sgn(M_X)=X$.  For an 
$\{+1,-1\}$-matrix $M$, the {\em sign rank} of $M$ is $\sr(M)=\mr(\sgn(M))$.
Paturi and Simon \cite{PS}, \cite[p. 106]{Lokam} showed that 
\[ \log_2\sr(M)\le \uppcc(M)\le  \log_2\sr(M)+1.\]  Thus the 
computation of $\sr(M)=\mr(\sgn(M))$ is of interest in the 
study of communication complexity.  A more thorough 
introduction to communication complexity and sign-rank  its connections to minimum rank 
are provided by Srinivasan's survey \cite{Sri} and Lokam's book 
\cite{Lokam}.

Forster \cite{F} established an important lower bound on the 
sign-rank of an $m\x n$ $\{+1,-1\}$-matrix.
\begin{thm}\label{Forst}{\rm \cite{F}} If $M$ is an $m\x n$ 
$\{+1,-1\}$-matrix, then
\[\sr(M)\ge  \frac   {\sqrt{mn}} {\| M \|},\]
where  $\| M\|$ is the spectral norm of $M$. 
\end{thm}
\noi An $n\x n$ Hadamard matrix $H$ realizes 
$ \sr(H)\ge  \frac   n {\sqrt n}=\sqrt n$ \cite{F}.

Some of the techniques described in Sections \ref{ssznz} and  
\ref{ssasymparam}  for nonzero patterns and loop digraphs 
(which are equivalent to square nonzero patterns) have been 
adapted to sign patterns. Triangle number used literally is 
less useful than the following generalization. An $n\x n$ sign 
pattern $X$ is {\em sign nonsingular} (SNS) if  every $n\x n$ 
real matrix
$B$ such that $\sgn(B)=X$ is nonsingular. The {\em SNS number} 
of a sign pattern $X$, denoted $\sns(X)$,  is the maximum size 
of an SNS sign pattern submatrix   of $X$ \cite{Hsp}. For a   
square sign pattern $X$,  the {\em (row) edit distance to 
nonsingularity}, $\ED(X)$, of $X$ is the minimum number of 
rows that must be changed to obtain an SNS pattern \cite{Hsp}.

\begin{thm}\label{snsED}{\rm \cite{Hsp}} For any $n\x n$ sign 
pattern $X$, $\sns(X)+ \ED(X)=n$.
\end{thm}

Sign patterns for which the minimum rank differs from the maximum 
rank by a fixed amount (such as 1) are discussed in  \cite{ AHLMG}.

\subsection{Trees}\label{sstrees} 
Trees were the first family of simple graphs for which the minimum 
rank problem was studied, and the minimum rank problem has been 
solved for square nonzero patterns and square sign patterns for 
which the graph (simple or loop, undirected or directed) of the 
nonzero positions is a tree.  
Minimum rank/maximum nullity can be computed by computing 
other parameters that are equal for trees.  
Since solving the minimum rank 
problem on connected components solves the problem, ``tree" can 
be replaced with  ``forest" throughout this discussion.

A {\em simple tree} is a connected acyclic simple graph.  
A {\em pseudocycle} is  a digraph  from which a cycle of 
length at least three can be obtained by reversing the 
direction of  
zero or more arcs.  A 
 {\em ditree} is a (simple or loop) digraph that does not 
 contain any pseudocycles.  A 
 {\em tree} is a graph that is 
 one of the following: a simple tree; a loop  
 graph that is a simple tree after all 
 loops are removed;  a ditree.  The 
 loop digraph $\G(X)$ of 
 an $n\x n$ sign pattern $X$ is equal to $\G(B)$ for 
 $B\in\Rnn$ such that $\sgn(B)=X$.  A square sign pattern 
 $X$ is a {\em tree sign pattern} if $\G(X)$ is a ditree.

It is well-known that that $\PC(T)=\M(T)$ for a 
simple tree $T$. In \cite{square1, cancun} the  definition 
of path cover number
is extended to 
graphs of other types 
and the analogous result established 
for trees of various types.  
In extending the definition of path cover, there is an issue of 
whether paths must be induced, 
which is irrelevant for trees, 
so here we 
extend the definition of path cover number only to  trees of various types.
 A loop (di)graph  $G$ {\em requires nonsingularity} if 
 $\M(G)=0$, i.e., $A\in\Q(G)$ implies $A$ is nonsingular 
 (this is analogous to sign nonsingularity); otherwise $G$ 
 allows singularity.  
Every simple graph allows singularity, which is immediate 
by considering $A-\lambda I$ where $A\in\Q(G)$ and  $\lambda$ 
is an eigenvalue of $A$.
In \cite[Definition 4.19]{square1}, the definition of path 
cover number was generalized to loop digraphs (and implicitly 
also to loop graphs) in a manner that retains the property 
$\PC(T)=\M(T)$ for a loop ditree.  A key idea was to ignore 
components that require nonsingularity (such components cannot 
exist in a simple graph). 
Let $T$ be a tree of any type. A {\em path cover}   of $T$ 
is a set of vertex disjoint  
paths whose deletion from $T$ leaves  a graph that requires 
nonsingularity (or the empty set).  The 
{\em path cover number} $\PC(T)$ is the minimum number of paths  
in a path cover.

\begin{thm}\label{MP}{\rm \cite{JL99, square1, cancun, Hsp}} 
For  a tree of any type or a tree sign pattern, $\M(T)=\PC(T)$.
\end{thm}

The parameters $\ZFS(T)$ and $\ED(T)$ are equal to $\M(T)$ when 
they have been defined. 

\begin{thm}\label{MZ}{\rm \cite{ AIM, square1, cancun}} For  a 
tree of any type, $\M(T)=\ZFS(T)$.
\end{thm}

\begin{thm}\label{mrtri}{\rm \cite{ square1}} For  loop 
ditree, $\M(T)=\ED(T)$ and $\mr(T)=\tri(T)$.
\end{thm}

\begin{thm}\label{ME}{\rm \cite{Hsp}} If $T$ is a 
tree sign pattern, $\M(T)=\ED(T)$ and $\mr(T)=\sns(T)$.
\end{thm}

For simple trees, the equality $\M(T)=\PC(T)$ was 
established  in \cite{JL99}, and was extended to 
$\M(T)=\PC(T)=\ZFS(T)$ in \cite{AIM}.   
The definition of   $\PC(T)$ was given for loop 
ditrees in \cite{square1}, where it was shown that 
a result in \cite{DHHHW}  implied $\M(T)=\PC(T)$ for 
loop trees, and  $\M(T)=\ZFS(T)=\ED(T)=\PC(T)$ was established 
for loop ditrees.  The equality $\M(T)=\ZFS(T)=\PC(T)$ was extended to 
simple ditrees in \cite{cancun} and for sign patterns $\M(T)=\ED(T)$ 
and $\mr(T)=\sns(T)$ were established by related methods in \cite{Hsp}.


\section{Related  problems described by (simple) graphs}\label{sother} 

\subsection{Minimum rank over other fields}\label{smrF}  

Recently there has been considerable interest in the study of 
minimum rank over fields other than the real numbers.    For a given graph $G$ of order $n$, let
\[\mrF(G)=\min\{\rank A: A\in\Fnn, A^T=A, \G(A)=G\}.\]

Graphs of minimum rank at most  2 over any field $F$ were characterized 
by a finite set of forbidden induced subgraphs in \cite{BHL04, BHL05}  
(with the set of forbidden subgraphs depending on the characteristic 
of $F$ and number of elements in $F$).
In \cite{DK} it was shown that  the set
of graphs of minimum rank at most $k$ over any finite field  is 
characterized by finitely many forbidden induced subgraphs.
In \cite{BGL} a complete set of forbidden induced subgraphs for minimum 
rank 3 over $\Z_2$ is determined.
In contrast to the finite field case, it is reported that an infinite 
set of forbidden induced subgraphs is needed to characterize minimum 
rank 3 over the real numbers \cite{Hall}.  Johnson, Loewy, and Smith 
characterize graphs having maximum nullity 2 over any infinite field 
\cite{JLS}.

In 2006 it was an open question whether the minimum rank over another 
field of characteristic zero (such as $\C$ or $\QQ$) could differ from 
$\mr(G)=\mr^\R(G)$  \cite{AIMrpt}.  In \cite{QRC} examples were given 
of graphs $G_1$ and $G_2$ such that $\mr^\R(G_1)>\mr^\C(G_1)$ and 
$\mr^\QQ(G_2)>\mr^\R(G_2)$.  Another example of a graph $G_3$ with 
$\mr^\QQ(G_3)>\mr^\R(G_3)$ was given in \cite{KR}.  The graphs $G_2$ 
and $G_3$ provided counterexamples to a conjecture in \cite{AHKLR}.

A {\em universally optimal matrix} is a (symmetric) integer matrix $A$ 
such that every off-diagonal entry of $A$ is 0, 1, or $-1$ (note for 
such a matrix $\G(A)$ is independent of field),  and for all fields 
$F$,  $\rank^F(A)=\mrF(\G(A))$ \cite{UOM}. 
In that paper universally optimal matrices 
were 
 used to show that a number of graphs in the the AIM Minimum 
Rank Graph Catalog \cite{AIMcat} have field independent minimum rank, 
and examples were presented to show that other graphs in the catalog 
are field dependent.  Additional results on universally optimal matrices 
and field independence are given in  \cite{Yeh2}.

\subsection{The graph parameter  $\eta(G)$}\label{seta}  

If $G$ is a  graph on vertices $\{1,\dots,n\}$, the {\em Haemers number} 
$\eta(G)$ is  defined to be the smallest
rank of any $n \x n$ matrix $B=[b_{ij}]$ (over any field) that 
satisfies $b_{ii}\ne 0$ for $i = 1, \dots, n$ and $b_{ij}=0$ if $i$ 
and $j$ are distinct nonadjacent vertices.  Clearly 
$\alpha(G)\le \eta(G)$ where $\alpha(G)$ is the independence number of $G$ (i.e., the maximum 
number of vertices with none adjacent).  The Laplacian matrix of $G$ 
shows that $\eta(G)\le n-c$ where $c$ is the number of connected components 
of order at least two.  Haemers has established a number of properties of 
$\eta(G)$, including   that $\eta(G)\le\chi({\overline G})$ 
(where $\chi(H)$ is the chromatic number of $H$), and $\eta(G)$ is an 
upper bound for the Shannon capacity of $G$ \cite{Hae81}. 

We now examine the relationship between $\eta(G)$ and the minimum rank 
parameters already discussed. Matrices satisfying the conditions of the 
Haemers number need not be symmetric but must have positive diagonal entries.
  If a symmetric matrix $A\in\Fnn$ satisfies the conditions of the 
  Haemers number for $G$, then $\G(A)$ is a subgraph of $G$.  The Haemers 
  number  $\eta(G)$ is not comparable to $\mr(G)$ as the next two examples 
  show.

\begin{ex}{\rm It is well known that $\mr(K_{1,3})=2$, and $\eta(K_{1,3})=3$
because $\alpha(K_{1,3})=3$.
}
\end{ex} 

\begin{ex}\label{exvcc}{\rm It is well known that $\mr(K_{3}\cp P_2)=3$, 
where $G\cp H$ denotes the Cartesian product (see \cite{AIM} for the definition).
If we number the vertices so that the two copies of $K_3$ are numbered $\{1,2,3\}$ and $\{4,5,6\}$, then we can see that $\eta(K_{3}\cp P_2)=2$ by considering the matrix 
$J_3\oplus J_3$ (where $J_3$ is the $3 \x 3$ matrix having every entry 
equal to 1).
}
\end{ex}

If $G$ is a connected graph, then any matrix $A\in\h_+(G)$ 
(see Section \ref{spsdmr}) satisfies 
the conditions on the matrices used to determine $\eta(G)$, so for a 
connected graph $G$, $\eta(G)\le \mrp^\C(G)$.  
A somewhat better upper bound for $\eta(G)$  is given by the asymmetric 
minimum rank of a loop digraph   (see Section 
\ref{ssasymparam}) obtained from $G$  
by replacing each edge by both arcs and adding a loop at each vertex, but this bound still requires a nonzero entry where 
an edge is present in the graph, and the Haemers number does not.  
Recall that the (edge) clique cover number $\cc(G)$ provides an upper 
bound for $\mrp(G)$.  The {\em vertex clique cover number}, i.e., the 
minimum number of cliques needed to cover all the vertices in $G$, is 
clearly an upper bound for $\eta(G)$; this was used in Example \ref{exvcc}.  
The vertex clique cover number can be much smaller than minimum rank. For 
example, $\mr(K_{n}\cp P_2)=n$ (and this does not change if asymmetric 
matrices are allowed), but the vertices of  $K_{n}\cp P_2$ can be covered 
by 2 cliques.

\subsection{Inverse inertia problem}\label{sinert} 

Barioli and Fallat \cite{BF} introduced the term {\em inertia balanced} to 
describe a graph with the property that there is 
a matrix that realizes the minimum rank of the graph and has
the number of negative eigenvalues equal to or one less than the number of 
positive eigenvalues.  Inertia balanced graphs played a crucial role in their 
study of the minimum rank of joins, and they showed that many graphs are 
inertia balanced.  They asked whether all graphs are inertia balanced.
Barrett, Hall, and Loewy \cite{BHL09} answered this question in the 
negative by exhibiting an example of a graph that is not inertia balanced.  

In  \cite{BHL09} they also began the study of the {\em inverse inertia problem},
i.e the question of determining what inertias are possible for matrices 
described by the graph.  For a given graph $G$, inverse inertia 
problem for $G$ lies in between  the minimum rank problem for $G$ and 
the inverse eigenvalue problem for $G$, i.e., the question of what 
spectra are possible for a matrix described by $G$.  Barrett, Hall, 
and Loewy \cite{BHL09} solved the inverse inertia problem for trees 
and provide a cut-vertex reduction formula for inverse inertia.  
The inverse inertia problem is solved for graphs of order at most 6 in 
  \cite{Betal}, where additional techniques for determining 
  inverse inertias are also presented.

\subsection{Minimum skew rank}\label{sskew} 

The majority of the work on minimum rank and related problems has 
focused on symmetric matrices.  There has also been work on 
matrices having a nonzero pattern 
described by a digraph, or having signs described by a sign pattern, 
see Section \ref{sasymmr}.  Recently there has also been interest in 
the problem of ranks of skew-symmetric matrices described by a graph. 
Such ranks are necessarily even, but full rank may not be possible 
(even in the case where the order of the graph is even).  
Let $\mrs(G)$ (respectively, $\Mrs(G)$) denote the minimum rank (maximum rank)
of matrices in the family $\sym_-(G)$ of real skew-symmetric matrices 
whose off-diagonal pattern of nonzero entries is described by the 
edges of $G$.  A {\em matching} is a set of edges such that all the 
vertices are distinct, $\match(G)$ denotes the number of edges in a 
maximum matching of $G$, and a matching is {\em perfect} if it 
includes every vertex.
\begin{thm}{\rm \cite{IMA}} Let $G$ be a graph.
\ben
\item Every even rank between $\mrs(G)$ and $\Mrs(G)$ can be realized.
\item $\mrs(G)=\ord G$ if and only if $G$ has a unique perfect matching.
\item $\Mrs(G)=2\match(G)$.
\item If $T$ is a tree, then $\mrs(T)=2\match(T)=\Mrs(T)$.
\item $\mrs(G)=2$ if and only if $G$ is a complete multipartite graph.
\een
\end{thm}
Minimum skew rank is computed for several families of graphs, the skew 
zero forcing number is defined, and related results over fields other 
than the real numbers are also presented in \cite{IMA}.


The  list above is intended to be used in conjunction with the 
bibliography in \cite{FH}, and references cited there are include 
here  only if they are cited in this paper.  \end{document}